\documentclass[12pt]{article}

\newcommand{\be}{\begin{equation}}
\newcommand{\ee}{\end{equation}}

\newtheorem{theorem}{Theorem}[section]
\newtheorem{lemma}[theorem]{Lemma}

\newtheorem{remark}[theorem]{Remark}

\begin{document}

\title{ Qualitative properties of solutions for system involving fractional Laplacian}
\author {{Ran Zhuo {\thanks{Partially supported by NSFC 11701207 and China Postdoctoral Science Foundation 2019M661473.}}
\hspace{.2in}Yingshu L\"{u}{\thanks{Partially supported by NSFC 11571233.}}}}

\date{}
\maketitle

\abstract{
In this paper, we consider the following nonlinear system involving the fractional Laplacian
\begin{equation}
\left\{\begin{array}{ll}
(-\Delta)^{s} u (x)= f(u,\,v),   \\
(-\Delta)^{s} v (x)= g(u,\,v),
\end{array}
\right. \label{ab1}
\end{equation}
in two different types of domains, one is bounded, and the other is unbounded, where $0<s<1$. To investigate the qualitative properties of solutions for fractional equations, the conventional methods are extension method and moving planes method. However, the above methods have technical limits in asymmetric and convex domains and so on. In this work, we employ
the direct sliding method for fractional Laplacian to derive the monotonicity of solutions for (\ref{ab1}) in $x_n$ variable in different types of domains.  Meanwhile, we develop a new iteration method for systems in the proofs which hopefully can be applied to solve other problems.}      


2010 MR Subject Classification: {35B50; 35J60; 35R11}      

\section{Introduction}
The work concerning qualitative properties of solutions was extensively investigated by many researchers.
Fruitful results have been obtained on the existence,
nonexistence, symmetry and regularity, and so on. Berestycki and Nirenberg \cite{BN} obtained the monotonicity for the unbounded positive solutions of elliptic equations in the case of a ``coercive`` Lipschitz graph. In \cite{BCN}, Berestycki, Caffarelli and Nirenberg proved the monotonicity and uniqueness for elliptic equations in unbounded Lipschitz domains. Angenent \cite{An} and  Cl\'{e}ment
and Sweers \cite{CS} derived that a bounded positive solution of elliptic equations only depends on the $x_n$-variable in a upper half space.
Chen, Li and Li \cite{CLL} worked on the symmetry and nonexistence of positive solutions of equations with fractional Laplacian in different types of domains. For more related results, please see \cite{CL, CL1, CFY, GNN, CT}, and the references therein.

In this work, we investigate the following system involving the fractional Laplacian:
\begin{equation}
\left\{\begin{array}{ll}
(-\Delta)^{s} u (x)= f(u,\,v),  \\
(-\Delta)^{s} v (x)= g(u,\,v)
\end{array}
\right. \label{in-1}
\end{equation}
in a bounded domain and an unbounded domain, respectively.

From an applicable point of view, the fractional Laplacian have caught researchers'
attention because of its nonlocality and its applications in physical sciences.
So far it has been utilized to model diverse physical phenomena, such as anomalous diffusion and quasi-geostrophic
flows, turbulence and water waves, molecular dynamics, and relativistic quantum mechanics
of stars (see \cite{BoG}, \cite{Co}, \cite{CV}, \cite{TZ} and the
references therein). It also has various applications in probability
and finance (\cite{A}, \cite{Be}). In particular, the fractional Laplacian
can be understood as the infinitesimal generator of a stable L\'{e}vy diffusion process \cite{Be}.

The fractional Laplacian is a nonlocal pseudo-differential
operator, taking the form
\begin{equation}
(-\Delta)^{s} u(x) = C_{n,\alpha} \, PV \int_{R^n}
\frac{u(x)-u(z)}{|x-z|^{n+2s}} dz,\label{Ad7}
\end{equation}
where $s\in(0,\,1)$ and \emph{PV} stands for the Cauchy principal value.
This operator is well defined for $u\in C^{\infty}_0(R^n)$. In this space,
it can also be defined equivalently in terms of the Fourier
transform
$$ (-\Delta)^{s} u (x) = {\cal{F}}^{-1}(|\xi|^{2s} {\cal{F}}u(\xi))(x), $$
where $\cal{F}$ is the Fourier transform, and ${\cal{F}}^{-1}$ is the inverse Fourier transform. The fractional Laplacian can be extended to locally integrable functions with certain growth control--the weighted $L^1$-space:
$$L_{2s}=\{u: R^n\rightarrow R \mid \int_{R^n}\frac{|u(x)|}{1+|x|^{n+2s}} \, d x <\infty\} \,\,\,\mbox{(see \cite{Si})}.$$

For $u\in L_{2s}$, we define $(-\Delta)^{s} u$ as a
distribution:
$$ ((-\Delta)^{s}u)( \phi) = \int_{R^n} u(x) (-\Delta)^{s} \phi(x)dx,  \;\;\; \forall \, \phi \in {C^{\infty}_c}(R^n).$$

To investigate the properties of solutions of equations involving the fractional Laplacian, the method of moving
planes in integral forms \cite{CL2} and extension method \cite{CS} have been general powerful tools. In \cite{ZCCY}, they employed the method of moving planes in integral forms (see \cite{CL}, \cite{CFY} and the references therein). This method requires the existence and some properties of Green's functions. But so far, there are few results about Green's functions in general domains. Meanwhile, it is much harder to prove the equivalence between the differential equations and the integral equations. In \cite{BCPS}, the authors used the Caffarelli-Silvestre's extention method. They technically required that $1/2<s<1$ and $u$ be globally bounded. In this paper, we employ the direct sliding method, which was developed by Berestycki and  Nirenberg \cite{BN}, to derive the monotonicity of solutions of (\ref{in-1}) in bounded and unbounded domains (also see \cite{WC}). By using the sliding method, we
weaken the global boundedness condition.

The analogue problem to (\ref{in-1}) for the fractional Laplacian has been investigated by many authors:
\begin{equation}
\left\{\begin{array}{ll}
(-\Delta)^{s} u_i (x)= f_i(u_1,\dots,u_m), \,\,\, x\in \Omega, \ & i=1,\dots,m, \\
u_i(x)=0, \,\,\,x\in {R}^n\setminus \Omega  \ &i=1,\dots,m.
\end{array}
\right. \label{ad-1}
\end{equation}
In the case $\Omega$ be a unit ball or half space, Mou \cite{M} proved the symmetry and monotonicity of positive solutions of (\ref{ad-1}) by the integral equation approach. When $i=1,2$, $f_1=u^p_2$, $f_2=u^q_1$, and $\Omega=R^n_+$, Quaas and Xia \cite{QX} obtained the non-existence
of positive solutions of (\ref{ad-1}) by the method of moving planes with an improved Aleksandrov-Bakelman-Pucci type estimate for the fractional Laplacian.

In this paper, we consider the nonlinear equations involving the fractional Laplacian in general domains. Due to the nonlocal nature of the operator, we need to set exterior conditions in domain $\Gamma$,
$$u(x)=\varphi(x),\,\,\,v(x)=\psi(x),\,\,\,x\in \Gamma^c.$$

In order to ensure the monotonicity of solutions, one has to impose the necessary exterior conditions {\bf{(P)}} on $\Gamma$:
 For any three points $x^1=(x',\,x^1_n)$, $x^2=(x',\,x^2_n)$ and $x^3=(x',\,x^3_n)$ with $x^1_n<x^2_n<x^3_n$, $x'\in {R}^{n-1}$,
and $x^1$, $x^3\in \Gamma^c$, $x^1$, $x^2$ and $x^3$ satisfy
$$\varphi(x^1)<u(x^2)<\varphi(x^3),\,\,\,\mbox{for }x^2\in \Gamma,$$
$$\varphi(x^1)\leq\varphi(x^2)\leq\varphi(x^3),\,\,\,\mbox{for }x^2\in \Gamma^c,$$
$$\psi(x^1)<v(x^2)<\psi(x^3),\,\,\,\mbox{for }x^2\in \Gamma,$$
$$\psi(x^1)\leq\psi(x^2)\leq\psi(x^3),\,\,\,\mbox{for }x^2\in \Gamma^c.$$

We first study (\ref{in-1}) in bounded domain $G$, and establish monotonicity of solutions for (\ref{in-1}) as following:
\begin{theorem}\label{thm1}
Let $u$, $v\in L_{2s}\cap C^{1,1}_{loc}(G)\cap C(\bar{G})$ and $(u,\,v)$ be a pair of solution for
\begin{equation}
\left\{\begin{array}{ll}
(-\Delta)^{s} u (x)= f(u,\,v), \,\,\,& x\in G,  \\
(-\Delta)^{s} v (x)= g(u,\,v), \,\,\,& x\in G,  \\
u(x)=\varphi(x),\,\,\,v(x)=\psi(x),\,\,\,& x\in G^c,
\end{array}
\right. \label{in-2}
\end{equation}
where $G$ is a convex bounded domain in $x_n$ direction, and $u$, $v$ satisfy the exterior condition {\bf{(P)}} on $G$. If $\frac{\partial f}{\partial v}>0$, $\frac{\partial g}{\partial u}>0$ for $x\in G$, and $\frac{\partial f}{\partial u}$, $\frac{\partial f}{\partial v}$, $\frac{\partial g}{\partial u}$, $\frac{\partial g}{\partial v}$ are bounded from above in $G$, then $u$ and $v$ are monotone increasing with respect to $x_n$-variable in $G$. More precisely, for any $\tau>0$, one has
$$u(x',\tau+x_n)>u(x',x_n),\,\,\, \forall (x',\tau+x_n), \,(x',x_n)\in G.$$
\end{theorem}

In the proof of Theorem \ref{thm1}, employing the argument by contradiction at extreme points,
we first derive the key tool-"narrow region principle in bounded domains". Combining with the sliding method, the solutions of system (\ref{in-2}) are proved to be monotonic in the bounded domain $G$. Next we briefly introduce the sliding method. For any positive real number $\tau$, by sliding downward $\tau$ units from the bounded domain $G$, we have
$$G_{\tau}=\{x-\tau e_n\,|\,x\in G\},$$
here $e_n=(0,\dots,0,1)$.   Denote
$$u^{\tau}(x)=u(x',x_n+\tau),\,\,\,v^{\tau}(x)=v(x',x_n+\tau),$$
$$\tilde{U}^{\tau}(x)=u^{\tau}(x)-u(x),\,\,\,\tilde{V}^{\tau}(x)=v^{\tau}(x)-v(x).$$

For $\tau$ sufficiently close to the width of $G$ in $x_n$ direction, it is easy to see that $G\cap G_{\tau}$
is a narrow region. Applying the "narrow region principle in bounded domains" yields that
\begin{equation}\tilde{U}^{\tau}(x)\geq 0,\,\,\,\tilde{V}^{\tau}(x)\geq 0,\,\,\,x\in G\cap G_{\tau}.\label{in-ad-1}
\end{equation}

Note that (\ref{in-ad-1}) provides a starting position to slide the domain $G_{\tau}$. Then we slide
$G_{\tau}$ back upward as long as inequality (\ref{in-ad-1}) holds to its
limiting position. In fact, the domain should be slid to $\tau=0$. We conclude that the solutions of (\ref{in-2}) are monotone increasing in $x_n$ variable.

Considering the unbounded domain $\Omega$, because of the unboundedness property for $\Omega$, the extremum points in $\Omega$ cannot be attained which makes it hard to apply the key tool-"narrow region principle in unbounded domains" directly. To overcome this difficulty, we first take a minimization sequence approaching to the infimum, and then make some perturbation about the sequence to attain
the extremum points in some bounded domain. Combining with the iteration method, we can deduce the "narrow region principle in unbounded domains". Then the monotonicity of solutions for (\ref{in-1}) in an unbounded domain $\Omega$ is obtained with the aid of the direct sliding method.

\begin{theorem}\label{thm2}
Let $u$, $v\in  L_{2s}\cap C^{1,1}_{loc}(R^n)$ be a pair of solution for
\begin{equation}
\left\{\begin{array}{ll}
(-\Delta)^{s} u (x)= f(u,\,v), \,\,\,& x\in \Omega,  \\
(-\Delta)^{s} v (x)= g(u,\,v), \,\,\,& x\in \Omega,  \\
u(x)=\varphi(x),\,\,\,v(x)=\psi(x),\,\,\,& x\in \Omega^c,
\end{array}
\right. \label{in-3}
\end{equation}
where $\Omega=\{x=(x',\,x_n)\in R^n\,|\,0<x_n<M\}$, $x'=(x_1,\,x_2,\dots,x_{n-1})$, and $M$ is a finite positive real number. $u,\,v$ satisfy the exterior condition {\bf{(P)}} on $\Omega$, and $u(x',\,\cdot)$ and $v(x',\,\cdot)$ are bounded with $x'\in R^{n-1}$. Suppose that $\frac{\partial f}{\partial v}>0$, $\frac{\partial g}{\partial u}>0$ for $x\in \Omega$, and $\frac{\partial f}{\partial u}$, $\frac{\partial f}{\partial v}$, $\frac{\partial g}{\partial u}$, $\frac{\partial g}{\partial v}$ are bounded from above in $\Omega$. Then $u(x)$ and $v(x)$ are monotone increasing in $x_n$-variable, that is, for any $\tau>0$,
$$u(x',\tau+x_n)>u(x',x_n),\,\,\, \forall (x',\tau+x_n), \,(x',x_n)\in \Omega.$$
\end{theorem}

\begin{remark}
Theorem \ref{thm2} still holds if $\Omega$ is any domain bounded in the $x_n$-direction.
\end{remark}

One of the interesting point about the monotonicity of solutions is that it helps to pave the
way for deriving existence, non-existence and useful Sobolev inequalities, as can be seen
in \cite{CL1}, \cite{CL2} and the references therein.

This paper is organized as follows. In section 2, we establish some lemmas, such as "narrow region principle in bounded domains" and so on.
In section 3 and 4, combining the lemmas in section 2 with the sliding method, we derive the monotonicity of solutions for (\ref{in-1}) on bounded domains and unbounded domains.

\section{Key tools in the sliding method}
The aim of this section is to show the key tools in the sliding method. More precisely, we investigate the narrow region principle in bounded domains and unbounded domains so that the sliding method can be initiated.

\begin{lemma} \label{lem1}
(Narrow region principle for system in bounded domains)
Let $\tilde{U}$, $\tilde{V}\in L_{2s}\cap C^{1,1}_{loc}(E)$ satisfy
\begin{equation}
\left\{\begin{array}{ll}
(-\Delta)^{s} \tilde{U} (x)-b_1(x)\tilde{U}(x)-c_1(x)\tilde{V}(x)\geq 0, \,\,\,& x\in E,  \\
(-\Delta)^{s} \tilde{V} (x)-b_2(x)\tilde{U}(x)-c_2(x)\tilde{V}(x)\geq 0, \,\,\,& x\in E,  \\
\tilde{U}(x)\geq 0,\,\,\,\tilde{V}(x)\geq 0,\,\,\,& x\in E^c,
\end{array}
\right. \label{lem1-bd}
\end{equation}
where $E$ is a bounded domain, $c_1(x)>0$, $b_2(x)>0$ in $E$, and $b_i$, $c_i$ are bounded from above in $E$, $i=1,\,2$.
Then for $d$ sufficiently small, which is the width of $E$ in $x_n$ direction, one has
\begin{equation}
\tilde{U}(x)\geq 0,\,\,\,\tilde{V}(x)\geq 0,\,\,\, x\in E.
\label{lem1-bd-1}
\end{equation}

\end{lemma}

{\bf{Proof of Lemma \ref{lem1}}}

If (\ref{lem1-bd-1}) is not valid, then at least one of $\tilde{U}$ and $\tilde{V}$ is less than zero at some point. We may assume that there exists a point $x^0\in E$ such that
$$\tilde{U}(x^0)=\min\limits_{R^n}\tilde{U}(x)<0.$$

By (\ref{lem1-bd}), for $d$ sufficiently small, we have
\begin{eqnarray}
0&\leq& (-\Delta)^{s} \tilde{U} (x^0)-b_1(x^0)\tilde{U}(x^0)-c_1(x^0)\tilde{V}(x^0) \nonumber\\
&\leq& \frac{c\tilde{U}(x^0)}{d^{2s}}-b_1(x^0)\tilde{U}(x^0)-c_1(x^0)\tilde{V}(x^0) \nonumber\\
&\leq& \frac{c\tilde{U}(x^0)}{d^{2s}}-c_1(x^0)\tilde{V}(x^0).\label{lem1-bd-2}
\end{eqnarray}

This implies that
\begin{equation}
\tilde{V}(x^0)<0, \label{lem1-bd-3}
\end{equation}
and
\begin{equation}
\tilde{V}(x^0)\leq \frac{c\tilde{U}(x^0)}{d^{2s}c_1(x^0)}. \label{lem1-bd-4}
\end{equation}

It follows from (\ref{lem1-bd-3}) that there exists some point $x^1\in E$ such that
$$\tilde{V}(x^1)=\min\limits_{R^n}\tilde{V}(x)<0.$$

Similar to (\ref{lem1-bd-2}), we derive that, for $d$ sufficiently small,
\begin{eqnarray}
0&\leq& (-\Delta)^{s} \tilde{V} (x^1)-b_2(x^1)\tilde{U}(x^1)-c_2(x^1)\tilde{V}(x^1) \nonumber\\
&\leq& \frac{c'\tilde{V}(x^1)}{d^{2s}}-b_2(x^1)\tilde{U}(x^0).\label{lem1-bd-5}
\end{eqnarray}
Moreover, we get
\begin{equation}
\tilde{U}(x^0)\leq \frac{c'\tilde{V}(x^1)}{d^{2s}b_2(x^1)}. \label{lem1-bd-6}
\end{equation}

Combining (\ref{lem1-bd-4}) and (\ref{lem1-bd-6}) yields that
\begin{equation}
\tilde{V}(x^0)\leq \frac{c c'\tilde{V}(x^0)}{d^{4s}c_1(x^0)b_2(x^1)}.
\label{lem1-bd-7}
\end{equation}
Thus, one has
\begin{equation}
1\geq\frac{c c'}{d^{4s}c_1(x^0)b_2(x^1)}.
\label{lem1-bd-8}
\end{equation}
(\ref{lem1-bd-8}) is impossible for sufficiently small $d$. Therefore (\ref{lem1-bd-1}) holds.

\begin{lemma}\label{lem2}
(Narrow region principle for system in unbounded domains)
Let $D_1=\{x=(x',x_n)\in R^n|0<x_n<2l\}$ be an unbounded narrow region with some bounded constant $l$, and $D_-=\{x=(x',x_n)\in R^n|x_n<0\}$. If ${U}$, ${V}\in L_{2s}\cap C^{1,1}_{loc}(D_1)$ satisfy
\begin{equation}
\left\{\begin{array}{ll}
(-\Delta)^{s} {U} (x)-\bar{b}_1(x){U}(x)-\bar{c}_1(x){V}(x)\geq 0, \,\,\,& x\in D_1,  \\
(-\Delta)^{s}{V} (x)-\bar{b}_2(x){U}(x)-\bar{c}_2(x){V}(x)\geq 0, \,\,\,& x\in D_1,  \\
{U}(x)\geq 0,\,\,\,{V}(x)\geq 0,\,\,\,& x\in D_-.
\end{array}
\right. \label{unbd}
\end{equation}
Suppose that $\bar{c}_1(x)>0$, $\bar{b}_2(x)>0$ in $D_1$, and $\bar{b}_i$, $\bar{c}_i$ are bounded from above in $D_1$, $i=1,\,2$.
Then for $l$ sufficiently small, we get
\begin{equation}
{U}(x)\geq 0,\,\,\,{V}(x)\geq 0,\,\,\, x\in D_1.
\label{unbd-1}
\end{equation}

\end{lemma}

{\bf{Proof of Lemma \ref{lem2}}}

 The argument, by contradiction, is standard. Suppose (\ref{unbd-1})
is false. Then at least one of ${U}(x)$ and ${V}(x)$ are less than zero at some points belonging to $D_1$. Without loss of generality,
we may assume that there are some points such that the values of $U$ at these points are less than zero. Then there exists a sequence $\{x^k\}^{\infty}_{k=1}\subset D_1$ such that
\begin{equation}
U(x^k)\rightarrow A=\inf\limits_{R^n}U(x)<0,
\label{unbd-2}
\end{equation}
with $|x^k_n|<l$, where $x^k_n$ is the $n-$th component of $x^k$.

Let
\begin{equation}
\eta(x)=\left\{\begin{array}{ll}
ae^{\frac{1}{|x|^2-l}}, \,\,\,& |x|<l,  \\
0, \,\,\,& |x|\geq l,
\end{array}
\right. \label{unbd-3}
\end{equation}
taking $a=e^{1/l}$ such that $\eta(0)=\max\limits_{R^n}\eta(x)=1$.

Set $\varphi_k(x)=\eta(x-x^k)$. Combining with (\ref{unbd-2}), there exists a positive sequence
$\{\epsilon^k\}^{\infty}_{k=1}$ such that
\begin{equation}
U(x^k)-\epsilon^k\varphi_k(x^k)<A<0,
\label{unbd-4}
\end{equation}
where $\epsilon^k\rightarrow 0$ as $k\rightarrow \infty$.

Obviously, for $x\in R^n\setminus B_l(x^k)$, $U(x)\geq A$ and $\varphi_k(x)=0$. Then we have
\begin{equation}
U(x^k)-\epsilon^k\varphi_k(x^k)< U(x)-\epsilon^k\varphi_k(x),\,\,\,\mbox{for }x\in R^n\setminus B_l(x^k),
\label{unbd-5}
\end{equation}
here $B_l(x^k)=\{x\in R^n\,|\,|x-x^k|<l\}$.

Define $U_k(x)=U(x)-\epsilon^k\varphi_k(x)$. It follows from (\ref{unbd-5}) that there exists some point $\bar{x}^k\in B_l(x^k)$ such that
$$U_k(\bar{x}^k)=\min\limits_{R^n}U_k(x)<0.$$

It is easy to see that
$$U(\bar{x}^k)\leq U(x^k),$$
and
$$U(\bar{x}^k)\rightarrow A,\,\,\,\mbox{as }k\rightarrow\infty.$$

Applying the first inequality of (\ref{unbd}) and the definition of the fractional Laplacian, we derive
\begin{eqnarray}
0&\leq& (-\Delta)^{s} {U} (\bar{x}^k)-\bar{b}_1(\bar{x}^k){U}(\bar{x}^k)-\bar{c}_1(\bar{x}^k){V}(\bar{x}^k) \nonumber\\
&=& (-\Delta)^{s} {U}_k (\bar{x}^k)-\bar{b}_1(\bar{x}^k){U}_k(\bar{x}^k)-\bar{c}_1(\bar{x}^k){V}(\bar{x}^k) \nonumber\\
&+& \epsilon^k (-\Delta)^{s}\varphi_k(\bar{x}^k)-\epsilon^k\bar{b}_1(\bar{x}^k)\varphi_k(\bar{x}^k) \nonumber\\
&\leq& (\frac{C}{l^{2s}}-\bar{b}_1(\bar{x}^k)) {U}_k(\bar{x}^k)-\bar{c}_1(\bar{x}^k){V}(\bar{x}^k)\nonumber\\
&+& \epsilon^k (-\Delta)^{s}\varphi_k(\bar{x}^k) -\epsilon^k \bar{b}_1(\bar{x}^k)\varphi_k(\bar{x}^k).
\label{unbd-6}
\end{eqnarray}
Then for sufficiently small $l$ and sufficiently large $k$, one has
\begin{eqnarray}
0\leq  \frac{C}{l^{2s}} {U}_k(\bar{x}^k)-\bar{c}_1(\bar{x}^k){V}(\bar{x}^k)+ o(\epsilon^k).
\label{unbd-7}
\end{eqnarray}

This implies that
\begin{equation}
{V}(\bar{x}^k)<0.\label{unbd-8}
\end{equation}

Based on (\ref{unbd-8}), there exists a sequence $\{z^k\}^{\infty}_{k=1} \subset D_1$ such that
\begin{equation}
V(z^k)\rightarrow B=\inf\limits_{R^n}V(x)<0,
\label{unbd-9}
\end{equation}

Set $\psi_k(x)=\eta(x-z^k)$. It is easy to see that
\begin{equation}
V(z^k)-\epsilon^k\psi_k(z^k)<B<0
\label{unbd-10}
\end{equation}
and
\begin{equation}
V(z^k)-\epsilon^k\psi_k(z^k)< V(x)-\epsilon^k\psi_k(x),\,\,\,\mbox{for }x\in R^n\setminus B_l(z^k),
\label{unbd-11}
\end{equation}

Define $V_k(x)=V(x)-\epsilon^k\psi_k(x)$. It follows that there exists some point $\bar{z}^k\in B_l(z^k)$ such that
$$V_k(\bar{z}^k)=\min\limits_{R^n}V_k(x)<0.$$

Similar to the proof of (\ref{unbd-6}), by the second inequality of (\ref{unbd}), we arrive at
\begin{eqnarray}
0 &\leq& (-\Delta)^{s}{V} (\bar{z}^k)-\bar{b}_2(\bar{z}^k){U}(\bar{z}^k)-\bar{c}_2(\bar{z}^k){V}(\bar{z}^k) \nonumber\\
&=& (-\Delta)^{s}{V}_k (\bar{z}^k)-\bar{b}_2(\bar{z}^k){U}_k(\bar{z}^k)-\bar{c}_2(\bar{z}^k){V}_k(\bar{z}^k) \nonumber\\
&-& \bar{b}_2(\bar{z}^k){\varphi}_k(\bar{z}^k) \epsilon^k + (-\Delta)^{s}{\psi}_k (\bar{z}^k) \epsilon^k- \bar{c}_2(\bar{z}^k){\psi}_k(\bar{z}^k) \epsilon^k \nonumber\\
&\leq & (\frac{C'}{l^{2s}}-\bar{c}_2(\bar{z}^k)){V}_k(\bar{z}^k) -\bar{b}_2(\bar{z}^k){U}_k(\bar{z}^k) \nonumber\\
&-& \bar{b}_2(\bar{z}^k){\varphi}_k(\bar{z}^k) \epsilon^k + (-\Delta)^{s}{\psi}_k (\bar{z}^k) \epsilon^k- \bar{c}_2(\bar{z}^k){\psi}_k(\bar{z}^k) \epsilon^k.
\label{unbd-12}
\end{eqnarray}
For sufficiently small $l$ and sufficiently large $k$, we derive
\begin{eqnarray}
0 \leq \frac{C'}{l^{2s}}{V}_k(\bar{z}^k) -\bar{b}_2(\bar{z}^k){U}_k(\bar{z}^k)+o(\epsilon^k).
\label{unbd-13}
\end{eqnarray}
That is
\begin{eqnarray}
{V}_k(\bar{z}^k)&\geq& C'l^{2s}\bar{b}_2(\bar{z}^k){U}_k(\bar{z}^k)+o(\epsilon^k) \nonumber\\
&\geq& C'l^{2s}\bar{b}_2(\bar{z}^k){U}_k(\bar{x}^k)+o(\epsilon^k).
\label{unbd-14}
\end{eqnarray}

By (\ref{unbd-7}), we have
\begin{eqnarray}
\bar{c}_1(\bar{x}^k){V}(\bar{x}^k)\leq \frac{C}{l^{2s}} {U}_k(\bar{x}^k)+ o(\epsilon^k).
\label{unbd-15}
\end{eqnarray}

Combining (\ref{unbd-14}) with (\ref{unbd-15}), we derive
\begin{eqnarray}
&&C'l^{2s}\bar{c}_1(\bar{x}^k)\bar{b}_2(\bar{z}^k){U}_k(\bar{x}^k) \nonumber\\
&\leq & \bar{c}_1(\bar{x}^k){V}(\bar{z}^k) \nonumber\\
&\leq & \bar{c}_1(\bar{x}^k){V}(\bar{x}^k) \nonumber\\
&\leq & \frac{C}{l^{2s}} {U}_k(\bar{x}^k)+ o(\epsilon^k).
\label{unbd-16}
\end{eqnarray}

This yields
\begin{eqnarray}
C'\bar{b}_2(\bar{z}^k)\bar{c}_1(\bar{x}^k)\geq \frac{C}{l^{4s}}.
\label{unbd-17}
\end{eqnarray}
For sufficiently small $l$, (\ref{unbd-17}) is impossible.

This completes the proof of Lemma \ref{lem2}.

\section{Monotonicity of solutions in bounded domains}

In this section, we will verify Theorem \ref{thm1}.

{\bf{The proof of Theorem \ref{thm1}}}

Consider the following system:
\begin{equation}
\left\{\begin{array}{ll}
(-\Delta)^{s} u (x)= f(u,\,v), \,\,\,& x\in G,  \\
(-\Delta)^{s} v (x)= g(u,\,v), \,\,\,& x\in G,  \\
u(x)=\varphi(x),\,\,\,v(x)=\psi(x),\,\,\,& x\in G^c,
\end{array}
\right. \label{bdmt-1}
\end{equation}
where $G$ is a convex bounded domain in $x_n$ direction, and we denote the width of $G$ in $x_n$ direction as $d$.

First we introduce some basic notations. For any positive real number $\tau$, denote
$$u^{\tau}(x)=u(x',x_n+\tau),\,\,\,v^{\tau}(x)=v(x',x_n+\tau),$$
$$G_{\tau}=\{x-\tau e_n\,|\,x\in G\},$$
here $x=(x_1,\dots,x_{n-1},x_n)=(x',x_n)$, $e_n=(0,\dots,0,1)$, and
$$\Sigma_{\tau}=G\cap G_{\tau}.$$

Define $$\tilde{U}^{\tau}(x)=u^{\tau}(x)-u(x),\,\,\,\tilde{V}^{\tau}(x)=v^{\tau}(x)-v(x).$$

This proof consists of two steps.

{\bf{Step 1. }}For $0<\tau<d$ sufficiently large, we want to show
\begin{equation}
\tilde{U}^{\tau}(x)\geq 0,\,\,\,\tilde{V}^{\tau}(x)\geq 0,\,\,\,\forall x\in G_{\tau}.
\label{bdmt-2}
\end{equation}
By the exterior conditions ({\bf{P}}) of $u$ and $v$, it is easy to see that
\begin{equation}
\tilde{U}^{\tau}(x)\geq 0,\,\,\,\tilde{V}^{\tau}(x)\geq 0,\,\,\,\forall x\in G_{\tau}\setminus\Sigma_{\tau}.
\label{bdmt-add-2}
\end{equation}

This implies that we only need to prove
\begin{equation}
\tilde{U}^{\tau}(x)\geq 0,\,\,\,\tilde{V}^{\tau}(x)\geq 0,\,\,\,\forall x\in\Sigma_{\tau}.
\label{bdmt-3}
\end{equation}

Applying the {\em{mean value theorem}} to the first equation of (\ref{bdmt-1}), we have
\begin{equation}
(-\Delta)^{s}\tilde{U}^{\tau}(x)= f_u(\xi^{\tau}_1,\,v^{\tau})\tilde{U}^{\tau}(x)+f_v(u,\,\zeta^{\tau}_1)\tilde{V}^{\tau}(x),\,\,\,x\in G,
\label{bdmt-4}
\end{equation}
where $\xi^{\tau}_1$ is between $u$ and $u^{\tau}$ in $G$, and $\zeta^{\tau}_1$ is between $v$ and $v^{\tau}$ in $G$.

Similarly, we have
\begin{equation}
(-\Delta)^{s}\tilde{V}^{\tau}(x)= g_u(\xi^{\tau}_2,\,v^{\tau})\tilde{U}^{\tau}(x)+g_v(u,\,\zeta^{\tau}_2)\tilde{V}^{\tau}(x),\,\,\,x\in G,
\label{bdmt-5}
\end{equation}
where $\xi^{\tau}_2$ is between $u$ and $u^{\tau}$ in $G$, and $\zeta^{\tau}_2$ is between $v$ and $v^{\tau}$ in $G$.

Note that $\Sigma_{\tau}$ is a narrow region for $0<\tau<d$ sufficiently large. Applying Lemma \ref{lem1} to $\tilde{U}^{\tau}$ and $\tilde{V}^{\tau}$ with
$$E=\Sigma_{\tau},\,\,\,b_1=f_u(\xi^{\tau}_1,\,v^{\tau}),\,\,\,c_1=f_v(u,\,\zeta^{\tau}_1),$$
$$b_2=g_u(\xi^{\tau}_2,\,v^{\tau}),\,\,\,c_2=g_v(u,\,\zeta^{\tau}_2),$$
we derive that (\ref{bdmt-3}) is valid. We conclude that (\ref{bdmt-2}) must hold.

{\bf{Step 2. }}Now we decrease $\tau$ as long as (\ref{bdmt-2}) holds to the limiting position. Define
$$\tau_0=\inf\{\tau\,|\,\tilde{U}^{\tau}(x)\geq 0,\,\tilde{V}^{\tau}(x)\geq 0,x\in\Sigma_{\tau},\,0<\tau<d\}.$$

We want to prove
\begin{equation}
\tau_0=0.\label{bd-0}
\end{equation}

If $\tau_0>0$, we can show that $G_{\tau_0}$ can be slid upward a little bit and we still have, for some small $\delta>0$ and $\tau\in(\tau_0-\delta,\,\tau_0)$,
\begin{equation}
\tilde{U}^{\tau}(x)\geq 0,\,\,\,\tilde{V}^{\tau}(x)\geq 0,\,\,\,\forall x\in \Sigma_{\tau}.
\label{bd-0'}
\end{equation}
This contradicts with the definition of $\tau_0$. Therefore (\ref{bd-0}) holds. We postpone proving (\ref{bd-0'}).

In fact, for $\tau_0>0$, we can show that
\begin{equation}
\tilde{U}^{\tau_0}(x)>0,\,\,\,\tilde{V}^{\tau_0}(x)>0,\,\,\, x\in \Sigma_{\tau_0}.
\label{bd-1'}
\end{equation}

Otherwise, at least one of $\min\limits_{x\in \Sigma_{\tau_0}}\tilde{U}^{\tau_0}(x)$ and $\min\limits_{x\in \Sigma_{\tau_0}}\tilde{V}^{\tau_0}(x)$ are equal to zero. We may assume that, there exists a point $\bar{x}\in \Sigma_{\tau_0}$ such that
$$\tilde{U}^{\tau_0}(\bar{x})=\min\limits_{x\in \Sigma_{\tau_0}}\tilde{U}^{\tau_0}(x)=0.$$

It follows from (\ref{bdmt-4}) that
\begin{eqnarray}
(-\Delta)^s \tilde{U}^{\tau_0}(\bar{x})=f_v(u(\bar{x}),\,\zeta^{\tau_0}_1(\bar{x}))\tilde{V}^{\tau_0}(\bar{x}).
\label{bd-1}
\end{eqnarray}

On the other hand, by the exterior condition ({\bf{P}}) of $u$, we arrive at
\begin{eqnarray}
(-\Delta)^s \tilde{U}^{\tau_0}(\bar{x})=C_{n,\,s}P.V.\int_{R^n}\frac{\tilde{U}^{\tau_0}(\bar{x})-\tilde{U}^{\tau_0}(y)}{|\bar{x}-y|^{n+2s}}dy<0.
\label{bd-2}
\end{eqnarray}

Combining (\ref{bd-1}) and (\ref{bd-2}) yields
$$\tilde{V}^{\tau_0}(\bar{x})<0.$$
This is a contradiction. Hence (\ref{bd-1'}) is valid. It follows that
\begin{equation}
\tilde{U}^{\tau_0}(x)>0,\,\,\,\tilde{V}^{\tau_0}(x)>0,\,\,\,\forall x\in \Sigma_{\tau_0}.
\label{bd-3}
\end{equation}

Next we can choose some closed $Q\subset \Sigma_{\tau_0}$ such that $\Sigma_{\tau_0}\setminus Q$ is a narrow region. Applying (\ref{bd-3}), we have
\begin{equation}
\tilde{U}^{\tau_0}(x)\geq c_0>0,\,\,\,\tilde{V}^{\tau_0}(x)\geq c_0>0,\,\,\,\forall x\in Q.
\label{bd-4}
\end{equation}

By the continuity of $\tilde{U}^{\tau}$ and $\tilde{V}^{\tau}$ in $\tau$, we obtain, for some small $\delta>0$ and $\tau\in(\tau_0-\delta,\,\tau_0)$,
\begin{equation}
\tilde{U}^{\tau}(x)\geq 0,\,\,\,\tilde{V}^{\tau}(x)\geq0,\,\,\,\forall x\in Q.
\label{bd-5}
\end{equation}

Applying the exterior condition ({\bf{P}}), we have, for some small $\delta>0$ and $\tau\in(\tau_0-\delta,\,\tau_0)$,
\begin{equation}
\tilde{U}^{\tau}(x)\geq 0,\,\,\,\tilde{V}^{\tau}(x)\geq0,\,\,\,\forall x\in \Sigma^c_{\tau}.
\label{bd-6}
\end{equation}
It follows from Lemma \ref{lem1} that for some small $\delta>0$ and $\tau\in(\tau_0-\delta,\,\tau_0)$
\begin{equation}
\tilde{U}^{\tau}(x)\geq 0,\,\,\,\tilde{V}^{\tau}(x)\geq0,\,\,\,\forall x\in \Sigma_{\tau}\setminus Q.
\label{bd-7}
\end{equation}

Combining (\ref{bd-5}), (\ref{bd-6}) and (\ref{bd-7}), we derive that, for some small $\delta>0$ and $\tau\in(\tau_0-\delta,\,\tau_0)$,
\begin{equation}
\tilde{U}^{\tau}(x)\geq 0,\,\,\,\tilde{V}^{\tau}(x)\geq0,\,\,\,\forall x\in \Sigma_{\tau}.
\label{bd-8}
\end{equation}
This implies (\ref{bd-0'}) holds. It follows that (\ref{bd-0}) must be true.

This completes the proof of Theorem \ref{thm1}.

\section{Monotonicity of solutions in unbounded domains}
In this section, we study system (\ref{in-3}). For convenience, we write down (\ref{in-3}) again:
\begin{equation}
\left\{\begin{array}{ll}
(-\Delta)^{s} u (x)= f(u,\,v), \,\,\,& x\in \Omega,  \\
(-\Delta)^{s} v (x)= g(u,\,v), \,\,\,& x\in \Omega,  \\
u(x)=\varphi(x),\,\,\,v(x)=\psi(x),\,\,\,& x\in \Omega^c,
\end{array}
\right. \label{unbdmt-1}
\end{equation}
where $\Omega=\{x=(x',\,x_n)\in R^n\,|\,0<x_n<M\}$, $x'=(x_1,\,x_2,\dots,x_{n-1})$. we will verify Theorem \ref{thm2}.

{\bf{Proof of Theorem \ref{thm2}.}}

First we introduce some necessary notations. For any $0\leq \tau \leq M$, set
$$u^{\tau}(x)=u(x',x_n+\tau),\,\,\,v^{\tau}(x)=v(x',x_n+\tau).$$
Let
$$\Omega_{\tau}=\{x-\tau e_n\,|\,x\in \Omega\},$$
which is obtained by sliding $\Omega$ downward $\tau$ units in $x_n$ direction, $e_n=(0,\,0,\dots,0,\,1)$.

Set $$D_{\tau}=\Omega\cap\Omega_{\tau},$$
$$U^{\tau}(x)=u^{\tau}(x)-u(x),\,\,\,V^{\tau}(x)=v^{\tau}(x)-v(x).$$

The proof consists of three steps.

{\bf{Step 1. }}For $0<\tau<M$ sufficiently large, we want to show that
\begin{equation}
U^{\tau}(x)\geq 0,\,\,\,V^{\tau}(x)\geq 0,\,\,\,x\in \Omega_{\tau}.
\label{un-1}
\end{equation}

Obviously,
$$\Omega_{\tau}=D_{\tau}\cup (\Omega_{\tau}\cap R^n_{-}).$$

By the exterior condition ({\bf{P}}) of $u$ and $v$, we get
\begin{equation}
U^{\tau}(x)\geq 0,\,\,\,V^{\tau}(x)\geq 0,\,\,\,x\in \Omega_{\tau}\cap R^n_{-}.
\label{un-2}
\end{equation}

It is easy to see that $u^{\tau}(x)$ and $v^{\tau}(x)$ satisfy the PDEs (\ref{unbdmt-1}). Combining with the {\em{mean value theorem}}, we obtain
\begin{equation}
(-\Delta)^{s}U^{\tau}(x)= f_u(\xi^{\tau}_1,\,v^{\tau})U^{\tau}(x)+f_v(u,\,\zeta^{\tau}_1)V^{\tau}(x),\,\,\,x\in \Omega,
\label{un-4}
\end{equation}
where $\xi^{\tau}_1$ is between $u$ and $u^{\tau}$ in $\Omega$, and $\zeta^{\tau}_1$ is between $v$ and $v^{\tau}$ in $\Omega$.

Similarly, we have
\begin{equation}
(-\Delta)^{s}V^{\tau}(x)= g_u(\xi^{\tau}_2,\,v^{\tau})U^{\tau}(x)+g_v(u,\,\zeta^{\tau}_2)V^{\tau}(x),\,\,\,x\in \Omega,
\label{un-5}
\end{equation}
where $\xi^{\tau}_2$ is between $u$ and $u^{\tau}$ in $\Omega$, and $\zeta^{\tau}_2$ is between $v$ and $v^{\tau}$ in $\Omega$.

For $\tau$ sufficiently close to $M$, $D_{\tau}$ is narrow region in $x_n$ direction. Applying the "narrow region principle for system on unbounded domains" (Lemma \ref{lem2}), we arrive at
\begin{equation}
U^{\tau}(x)\geq 0,\,\,\,V^{\tau}(x)\geq 0,\,\,\,\forall x\in D_{\tau}.
\label{un-6}
\end{equation}

Combining (\ref{un-2}) and (\ref{un-6}), we derive that (\ref{un-1}) must hold.

{\bf{Step 2. }}(\ref{un-1}) provides a starting point to carry out the sliding method. Now we decrease $\tau$ as long as (\ref{un-1}) holds to the
limiting position. Define
$$\tau_0=\inf\{\tau\,|\,U^{\tau}(x)\geq 0,\,V^{\tau}(x)\geq 0,\,x\in D_{\tau},\,0<\tau<M\,\}.$$

We will show that
\begin{equation}
\tau_0=0.\label{un-7}
\end{equation}

Otherwise, suppose that $\tau_0>0$, we can show that $\Omega_{\tau}$ can be slid upward a little bit and we still have, for some small $\delta>0$,
\begin{equation}
U^{\tau}(x)\geq 0,\,\,\,V^{\tau}(x)\geq 0,\,\,\,\tau_0-\delta<\tau\leq \tau_0.\label{un-8}
\end{equation}
This is a contradiction with the definition of $\tau_0$. Then (\ref{un-7}) holds. We delay to proof (\ref{un-8}).

To prove (\ref{un-8}), we first show that
\begin{equation}
\inf_{x\in D_{\tau_0}}U^{\tau_0}(x)>0,\,\,\,\inf_{x\in D_{\tau_0}}V^{\tau_0}(x)>0.\label{un-9}
\end{equation}

If (\ref{un-9}) is not true, then at least one of $\inf_{x\in D_{\tau_0}}U^{\tau_0}(x)$ and $\inf_{x\in D_{\tau_0}}V^{\tau_0}(x)$ is
equal to zero. We may assume that
$$\inf_{x\in D_{\tau_0}}U^{\tau_0}(x)=0.$$

Hence, there exists a sequence $\{x^k\}^{\infty}_{k=1}\subset D_{\tau_0}$ such that
\begin{equation}
U^{\tau_0}(x^k)\rightarrow 0,\,\,\,\mbox{as }k\rightarrow \infty.
\label{un-10}
\end{equation}

Set
\begin{equation}
\eta(x)=\left\{\begin{array}{ll}
ae^{\frac{1}{|x|^2-r}}, \,\,\,& |x|<r,  \\
0, \,\,\,& |x|\geq r,
\end{array}
\right. \label{un-12}
\end{equation}
choosing $a=e^{1/r}$ and $r=\frac{M-\tau_0}{2}$, such that $\eta(0)=\max\limits_{R^n}\eta(x)=1$.

Let $\varphi_k=\eta(x-x^k)$. There exists a positive sequence $\{\varepsilon_k\}$ such that
$$U^{\tau_0}(x^k)-\varepsilon_k \varphi_k(x^k)<0$$
with $\varepsilon_k\rightarrow 0,$ as $k\rightarrow\infty$.

For any $x\in D_{\tau_0}\setminus B_r(x^k)$, $U^{\tau_0}(x)\geq 0$ and $\varphi_k(x)=0$. It is easy to see that,
\begin{equation}
U^{\tau_0}(x^k)-\varepsilon^k\varphi_k(x^k)< U^{\tau_0}(x)-\varepsilon^k\varphi_k(x),\,\,\,\mbox{for }x\in D_{\tau_0}\setminus B_r(x^k),
\end{equation}
where $B_r(x^k)=\{x\in R^n\,|\,|x-x^k|<r\}$.

It follows that there exists some point $\tilde{x}^k\in B_r(x^k)\cap D_{\tau_0}$ such that
\begin{equation}
U^{\tau_0}(\tilde{x}^k)-\varepsilon^k\varphi_k(\tilde{x}^k)=\min_{x\in D_{\tau_0}}(U^{\tau_0}(x)-\varepsilon^k\varphi_k(x))<0.
\label{un-13}
\end{equation}

Combining (\ref{un-10}) and (\ref{un-13}) yields that
$$U^{\tau_0}(x^k)\leq U^{\tau_0}(\tilde{x}^k)\leq U^{\tau_0}(x^k)-\varepsilon^k\varphi_k({x}^k)+\varepsilon^k\varphi_k(\tilde{x}^k).$$
Obviously, as $k\rightarrow \infty$,
\begin{equation}
U^{\tau_0}(\tilde{x}^k)\rightarrow 0.\label{un-14}
\end{equation}

By (\ref{un-4}), we derive that, for $k$ sufficiently large,
\begin{eqnarray}
&&(-\Delta)^s (U^{\tau_0}-\varepsilon_k\varphi_k)(\tilde{x}^k)  \nonumber\\
&=& f_u(\xi^{\tau_0}_1(\tilde{x}^k),\,v^{\tau_0}(\tilde{x}^k))U^{\tau_0}(\tilde{x}^k)
+f_v(u(\tilde{x}^k),\,\zeta^{\tau_0}_1(\tilde{x}^k))V^{\tau_0}(\tilde{x}^k)+ o(\varepsilon_k),\label{un-16}
\end{eqnarray}
where $\xi^{\tau_0}_1(\tilde{x}^k)$ is between $u(\tilde{x}^k)$ and $u^{\tau_0}(\tilde{x}^k)$, and $\zeta^{\tau_0}_1(\tilde{x}^k)$ is between $v(\tilde{x}^k)$ and $v^{\tau_0}(\tilde{x}^k)$.

On the other hand, employing the definition of the fractional Laplacian,
\begin{eqnarray}
&&(-\Delta)^s (U^{\tau_0}-\varepsilon_k\varphi_k)(\tilde{x}^k)  \nonumber\\
&=& C_{n,s} P.V.\int_{R^n} \frac{(U^{\tau_0}-\varepsilon_k \varphi_k)(\tilde{x}^k)-(U^{\tau_0}-\varepsilon_k \varphi_k)(y)}{|\tilde{x}^k-y|^{n+2s}}dy \nonumber\\
&\leq & c \int_{B^c_r(\tilde{x}^k)} \frac{(U^{\tau_0}-\varepsilon_k \varphi_k)(\tilde{x}^k)-(U^{\tau_0}-\varepsilon_k \varphi_k)(y)}{|\tilde{x}^k-y|^{n+2s}}dy \nonumber\\
&\leq & c \int_{B^c_r(0)} \frac{(U^{\tau_0}-\varepsilon_k\varphi_k)(\tilde{x}^k)-(U^{\tau_0}-\varepsilon_k\varphi_k)(y+\tilde{x}^k)}{|y|^{n+2s}}dy \nonumber\\
&\leq & c \int_{B^c_r(0)} \frac{-U^{\tau_0}(y+\tilde{x}^k)}{|y|^{n+2s}}dy. \label{un-17}
\end{eqnarray}

Set $u_k(x)=u(x+\tilde{x}^k)$, $U^{\tau}_k(x)=U^{\tau}(x+\tilde{x}^k)$. By Arzel$\grave{a}$-Ascoli theorem, we have
$$u_k(x)\rightarrow u_{\infty}(x),\,\,\,\mbox{as }k\rightarrow \infty,\,\,\,\mbox{in }R^n.$$

Hence, as $k\rightarrow \infty$,
\begin{equation}
U^{\tau}_k(x)\rightarrow U^{\tau}_\infty(x)=u^{\tau}_{\infty}(x)-u_{\infty}(x),\,\,\,x\in B^c_r(0).\label{un-18}
\end{equation}

Combining (\ref{un-16}), (\ref{un-17}) and (\ref{un-18}), we deduce that, as $k\rightarrow \infty$,
\begin{equation}
\int_{B^c_r(0)} \frac{-U^{\tau_0}_{\infty}(y)}{|y|^{n+2s}}dy \geq 0.\label{un-19}
\end{equation}

Obviously, (\ref{un-19}) holds unless
\begin{equation}
U^{\tau_0}_{\infty}(y)\equiv 0,\,\,\,y\in B^c_r(0).
\label{un-20}
\end{equation}

By (\ref{un-20}), we derive
\begin{equation}
u_{\infty}(x',x_n)=u_{\infty}(x',x_n+\tau_0)=\dots=u_{\infty}(x',x_n+m\tau_0)
\label{un-22}
\end{equation}
for any $m\in {N}^+$.

Choosing $(x',x_n)\in \Omega$, and taking $m$ large enough such that $(x',x_n+m\tau_0)\in \Omega^c$, we apply the exterior condition on $u$
to derive a contradiction with (\ref{un-22}). Thus (\ref{un-9}) holds.

Choosing sufficiently large $K\subset D_{\tau_0}$ such that $D_{\tau_0}\setminus K$ is narrow in $x_n$ direction. Combining (\ref{un-9})
with the continuity of $U^{\tau}$ and $V^{\tau}$ in $\tau$, we derive that, for some small $\delta>0$,
\begin{equation}
U^{\tau_0-\delta}(x)\geq 0,\,\,\,V^{\tau_0-\delta}(x)\geq 0,\,\,\,x\in K.
\label{un-23}
\end{equation}

Meanwhile, applying the exterior condition ({\bf{P}}), we have
\begin{equation}
U^{\tau_0-\delta}(x)\geq 0,\,\,\,V^{\tau_0-\delta}(x)\geq 0,\,\,\,x\in (D_{\tau_0-\delta})^{c}.
\end{equation}

Employing Lemma \ref{lem2}, we derive
\begin{equation}
U^{\tau_0-\delta}(x)\geq 0,\,\,\,V^{\tau_0-\delta}(x)\geq 0,\,\,\,x\in D_{\tau_0-\delta} \setminus K.
\label{un-24}
\end{equation}

Combining (\ref{un-23}) with (\ref{un-24}), we obtain
\begin{equation}
U^{\tau_0-\delta}(x)\geq 0,\,\,\,V^{\tau_0-\delta}(x)\geq 0,\,\,\,x\in D_{\tau_0-\delta}.
\end{equation}

This contradicts the definition of $\tau_0$. Hence (\ref{un-7}) is valid. We conclude that $u$ and $v$ are increasing in $x_n$ variable.

This completes the proof of Theorem \ref{thm2}.

\end{document}